\documentclass[12pt]{article}
\usepackage{amsmath,amssymb,amsthm,amsfonts}
\usepackage{latexsym}
\begin{document}

\newtheorem{prop}{Proposition}[section]
\newtheorem{cor}{Corollary}[section] 
\newtheorem{theo}{Theorem}[section]
\newtheorem{lem}{Lemma}[section]
\newtheorem{rem}{Remark}[section]
\newtheorem{con}{Conjecture}[section]

\renewcommand{\theequation}{\thesection.\arabic{equation}} 
\setcounter{page}{1} 
\noindent 
\begin{center}
{\Large \bf Merging percolation on $Z^d$ and classical random graphs:
Phase transition}
\end{center}

\begin{center}    
TATYANA S. TUROVA\footnote{Research was 
supported by the Swedish Natural Science Research
Council} and THOMAS VALLIER
\end{center}

\begin{center}
{\it Mathematical Center, University of
Lund, Box 118, Lund S-221 00, 
Sweden. }
\end{center}

\begin{abstract}
We study a  random graph model which is a superposition of
the bond percolation model on $Z^d$ with probability $p$ of an edge, and a classical random graph $G(n,
c/n)$. We show that
this model, being a {\it homogeneous} random
graph,
has a natural relation to the so-called
"rank 1 case" of {\it inhomogeneous} random
graphs. This allows us to use the newly developed theory of
inhomogeneous random graphs to describe the phase diagram on the set of parameters $c\geq 0$ and $0 \leq p<p_c$, where $p_c=p_c(d)$ is the critical probability for the bond percolation on $Z^d$.
The phase transition is similar to the classical random graph, it is
of the second order. We also find the scaled size of the largest connected
component above the phase transition.

\end{abstract}

\section{Introduction.}
\setcounter{equation}{0}

We consider a graph on the set of vertices $B(N):=\{-N, \ldots, N\}^d$
in $Z^d$, $d\geq 1$,
with two types of edges: the short-range edges connect independently 
 with probability $p$ each pair $u$ and $v$  if 
$|u-v|=1$, and the long-range edges
 connect independently any pair of two vertices with
probability $c/|B(N)|$. (Here 
for any set $A$ we denote $|A|$ the number of the elements in $A$.)
This graph, call it $G_N(p,c)$ is a superposition of the bond percolation model (see, e.g., \cite{G}),
where
each pair of
neighbours in $Z^d$ is connected with probability $p$, and a random
graph model $G_{n, c/n}$ (see, e.g., \cite{JLR}) on $n$ vertices,
where each vertex is connected to any other vertex with
probability $c/n$; all the edges in both models are independent. By this definition there can be one or two  edges between two vertices in graph $G_N(p,c)$, and in the last case the edges are of different types.

The introduced model is a simplification of the most
common graphs designed to study natural  phenomena, in particular,  biological 
neural networks \cite{TV}.
Notice the difference between  $G_N(p,c)$ and the so-called
"small-world" models intensively studied after \cite{WS}. In the 
 "small-world" models where edges from the grid  may be kept or
 removed,  only finite number (often at most  $2d$) of the long-range edges may
 come out of each vertex, and the probability of those is a  fixed
 number.
 
We are interested in the connectivity of the introduced graph $G_N(p,c)$ as $N
\rightarrow \infty$. We say that two vertices are connected, if there is a path of edges, no matter of which types, between them.
Clearly, if $c=0$, we have a purely bond percolation model on $Z^d$, where any edge from the grid is kept (i.e., "is open" in the terminology of percolation theory) with a probability $p$, or, alternatively, removed with a probability $1-p$. 
Let us recall some basic facts from the percolation theory which we need here.
Denote $C$ an open cluster containing the origin of $Z^d$ in the bond percolation model. It is known (see, e.g., \cite{G})
 that for any $d \geq 1$ there is $p_c =p_c(d)$ such that 
\begin{displaymath}
    {\bf P}\{|C|=\infty\} \ \left\{ \begin{array}{ll}
            =0, & \mbox{ if } p<p_c, \\ \\
            >0, & \mbox{ if } p>p_c ,
        \end{array} \right.
\end{displaymath}
where $0<p_c<1$, unless $d=1$, in which case, obviously, $p_c=1$.
We shall assume here that $0<p<p_c$, in which case the connected components
formed by the short-range edges only, are finite with probability one.
Recall also that for all $0<p<p_c$
the limit 
\begin{equation}\label{X}
\zeta(p)=lim_{n \rightarrow \infty}\left(
-\frac{1}{n} \log {\bf P} \left\{|C|=n \right\}\right) 
\end{equation}
exists and satisfies $\zeta(p)>0$ (Theorem (6.78) from \cite{G}).

Let further $ C_1 \Big( G \Big)$ denote the size (the number of vertices) of the largest
connected component in a graph $G$. 

\begin{theo}\label{T}
Assume, that  $d\geq 1$ and  $0\leq p < p_c(d)$.
Define
\begin{equation}\label{cr1}
  c^{cr} (p) = \frac{1}{{\bf E}|C|}.
\end{equation}

\noindent
$\mbox{i)}$ If 
$c  < c^{cr} (p)$ 
define $y$ to be the root of
${\bf E} \,c|C| e^{c|C|y }=1$, and set
\begin{equation}\label{alp}
  \alpha(p,c):=\left(c+cy-{\bf E} \, ce^{c|C|y}\right)^{-1}.
\end{equation}
Then for any $\alpha > \alpha(p,c)$
\begin{equation}\label{E1}
 {\bf P} \left\{ C_1  \Big( G_N(p,c)  \Big) >  
\alpha \log |B(N)|
\right\}  \rightarrow  0 .
\end{equation}
 as  $N
\rightarrow \infty$.

\medskip

\noindent
$\mbox{ii)}$ If 
$c  \geq  c^{cr} (p)$ then 
\begin{equation}\label{lp}
\frac{C_1 \Big(
G_N(p,c)
    \Big)}{|B(N)| }
\stackrel{P}{\rightarrow} \beta  
\end{equation}
as $ N \rightarrow \infty$,
with $\beta = \beta (p,c)$
defined as the maximal solution to
\begin{equation}\label{be}
  \beta = 1- {\bf E} \left\{ 
      e^{-{c}\, \beta \, |C|}
  \right\}.
\end{equation}
\end{theo}

\bigskip

In view of (\ref{X}) it is obvious
 that $\chi(p):={\bf E}|C|<\infty$ for all 
$0\leq p < p_c$. It is also known (see Theorem (6.108) and (6.52) in \cite{G}) 
that $\chi(p)$ is analytic function of $p$ on $[0,p_c)$ and $\chi(p) \rightarrow \infty$ as $p \rightarrow p_c$. 
This implies that $c^{cr}(p)$ is continuous, strictly decreasing function on $[0,p_c)$ with
$c^{cr}(0)=1$ and $c^{cr}(p_c)= 0$. 
Hence, $c^{cr}$ has inverse, i.e., for any $0<c<1$ there is a unique
$0<p^{cr}(c)<p_c=p_c(d)$
such that $c^{cr}\left(p^{cr}(c)\right)=c$. 
This leads to the following duality of the result in Theorem \ref{T}.

\begin{cor}
For any $0<c<1$ there is a unique
$0<p^{cr}(c)<p_c$ such that
for any $p^{cr}(c)<p<p_c$ graph $G_N(p,c)$ has a giant
component with a size $O(|B(N)|)$  {\bf whp} (i.e., with probability tending to one
as $N\rightarrow \infty$),
but for any  $p<p^{cr}(c)$ the size of the largest connected
component in  $G_N(p,c)$  {\bf whp} is at most  $O(\log |B(N)|)$.\hfill$\Box$
  \end{cor}
  
 Hence, 
Theorem \ref{T} may also tell us something  about the "distances"
between the components of a random graph when it is considered on the
vertices of ${\bf Z}^d$.

It is worth mentioning that
the symmetry between $c^{cr}$ and $p^{cr}$ is most spectacular in the
dimension one  case, when $p_c(1)=1$.
Notice, that $d=1$ case is exactly solvable,
and this is the only case when we know the distribution of
$|C|$:
\begin{equation}\label{Cd1}
{\bf P} \left\{|C|=k \right\}=(1-p)^2kp^{k-1},  \ \ \ k \geq 1.
\end{equation}
Hence, if $d=1$ we compute for all $0\leq p < 1=p_c(1)$
\begin{equation}\label{crd1}
c^{cr} (p) = \frac{1-p}{1+p},
\end{equation}
which also yields
\[p^{cr} (c) = \frac{1-c}{1+c},\]
for all $0\leq c < 1$. (For more details on  $d=1$ case we refer to
\cite{TT}.)

\begin{rem}
For any fixed $c$ function $\beta (p,c)$ is continuous at $p=0$:
 if  $p=0$, i.e., when our graph is merely a classical $G_{n,c/n}$
  random graph, then $|C| \equiv 1$ and  (\ref{be}) becomes a
  well-known relation.
\end{rem}

Furthermore,
for any fixed $c$, if $p=0$, it is not difficult to derive from
(\ref{alp})
that 
\begin{equation}\label{aoc}
\alpha (0,c) = \frac{1}{c-1+|\log c|}.
\end{equation}
But $\log n/(c-1+|\log c|)$
is known (see Theorem 7a in \cite{ER}) to be the principal term in the asymptotics (in probability) 
of the largest connected component of $G_{n,c/n}$.
This inevitably leads to the following (open) question. Will 
\begin{equation}\label{E2}
 {\bf P} \left\{ C_1  \Big( G_N(p,c)  \Big) < 
\alpha_1 \log |B(N)|
\right\}  \rightarrow  0 
\end{equation}
 hold also for all  $\alpha_1 <
\alpha(p,c)$ when $0<p<p_c$ and $c<c^{cr}(p)$?

  \bigskip
   It is easy to check that if $c \leq c^{cr}$ then the equation \eqref{be} does not have a
  strictly  positive
  solution, while $\beta=0$ is always a solution to \eqref{be}. 
This allows us to  derive
  \begin{equation}\label{be'}
  \beta '_{c} \mid _{\ c \downarrow c^{cr} }=2\frac{\left({\bf E} |C|\right)^{3}}{ {\bf E} (|C|^2)} >0,
\end{equation}
which confirms  that the phase
transition remains to be  of the second order  for any $p<p_c$, as it is for
$p=0$, i.e., in the case of classical random graph.

More similarities and differences between  our model and the "mean-field"
case one can see in  the following example. Let $d=1$, in which case
(\ref{Cd1})
holds. Introducing a random variable $X$ with the first-success distribution
\[{\bf P} \left\{ X=k \right\} =(1-p)p^{k-1}, \ \ \ \ k=1,2,
  \ldots \ ,\]
one can rewrite (\ref{be}) as follows (see the details in \cite{TT})
\[
\beta = 1-\frac{1}{{\bf E} X}\ {\bf E} \left\{ X
      e^{-{c}X\ \beta}
  \right\}.
\]
This equation looks 
somewhat  similar to the equation obtained in \cite{CL} for the "volume" (the sum of
degrees of the involved vertices)
of the giant component in the graph with a given sequence of the
expected degrees.
Note, however, that  in our model
 the critical mean degree 
when $c=c^{cr}$ and $N\rightarrow \infty$   is given
according to (\ref{crd1}) by
\begin{equation}\label{deg}
2p+c^{cr} = 2p + \frac{1-p}{1+p} =
1 +  \frac{2p^2}{1+p} 
  \end{equation}
which  is strictly greater than 1 for all positive $p<1$.
This is in a contrast with the model studied in \cite{CL},
where the
critical expected average degree is still 1 as  in the classical random graph.

Although our model (it can be considered on a torus, in the limit the result is the same)
is a perfectly {\it homogeneous} random graph, in the
sense that the degree distribution is the same for any  vertex, we
study it via {\it   inhomogeneous } random graphs, making use of the recently
developed theory from \cite{BJR}. The idea is the following. First, we consider  the
subgraph induced by the short-range edges, i.e., the edges which connect two neighbouring nodes with probability $p$. It is composed of the
 connected clusters (which may consist  just of one single vertex) in
$B(N)$. Call a
{\it macro-vertex} each of the connected components of this subgraph.
 We say that
a macro-vertex is of type $k$, if $k$ is the number of
vertices in it. 
Conditionally on the set of  macro-vertices, we
consider a graph on these macro-vertices induced by the long-range
connections. Two macro-vertices are said to be connected if there is
at least one (long-range type) edge  between two vertices
belonging to different macro-vertices.  
Thus the probability of an edge between two macro-vertices $v_i$ and $v_j$
of types $x$
and $y$ correspondingly,  is 
\begin{equation}\label{p}
    {\widetilde p}_{xy}(N):= 1-\left(1-\frac{c}{|B(N)|} \right)^{x y}.
\end{equation}
Below we argue that this model fits the conditions of a general
inhomogeneous graph model defined in \cite {BJR},
find the critical parameters
and characteristics for the
graph on macro-vertices,  and then we turn back to the original model. 
We use essentially the results from \cite {BJR} to derive
(\ref{be}). The 
result in the subcritical phase (part $i)$ of Theorem \ref{T}) does
not follow by the theory in \cite {BJR}; we discuss this in the end of
Section 2.4.

Notice also that
in order to analyze the introduced model, we 
derive here some result on the joint distribution of the sizes of
clusters in the percolation model (see Lemma \ref{L1} below), which
 may be of interest on its own.

The principle of treating  some local structures in a graph as new
vertices ("macro-vertices"),  and then
considering a graph induced by the original model on these vertices
appears to be rather general. For example,  in \cite{J}
a different graph model was also  put into a framework of
inhomogeneous graphs theory by
certain restructuring. This method should be useful for
 analysis of a broad class of complex structures, whenever  one can
 identify local
and global connections.  Some examples of such models one can find in
\cite{M}.

Finally we comment that
our result should help to study
a  model for the propagation of the neuronal activity
introduced in \cite{TV}. Here we show that
 a giant component in the graph can emerge from two sources, none of
 which can be neglected, but each of
 which may be in the subcritical phase, i.e., even  when both 
 $p<p_c$ and
 $c<1$. In particular, for any $0<c<1$ we can find $p<p_c$ which allows
 with a positive probability the
 propagation of  impulses  through the large part of the network due
 to the local activity.

\section{Proof}
\setcounter{equation}{0}

\subsection{Random graph on macro-vertices.}

  Consider  now the
subgraph on $B(N)$ induced by the short-range
edges only, which is a purely bond percolation model.
 By the
construction
this subgraph, call it $G_N^{(s)}(p)$,
is composed of
a random number of
 clusters (of connected vertices) of random  sizes. We call here the size of a cluster the
number of its vertices (it may be just one). We recall here more results from
percolation theory which we shall use later on.

Let $K_N$ denote the number of the connected components (clusters) in
$G_N^{(s)}(p)$,
and let 
\begin{equation}\label{t} 
{\bf X} = \{X_1, X_2, \ldots , X_{K_N}\}
\end{equation}
denote the collection  of all  connected clusters $X_i$ in
$G_N^{(s)}(p)$.
We shall also use $X_i$ to denote the set of vertices in the
$i-$th cluster.
By this definition  $\sum _{i=1}^{K_N}|X_i|=|B(N)|$.

\medskip

\noindent
{\bf Theorem} [\cite{G}, (4.2) Theorem, p. 77]{\it
\begin{equation}\label{P1}
\frac{K_N}{|B(N)|} \ {\rightarrow} \ \kappa(p):= {\bf
  E}\frac{1}{|C|}
\end{equation}
a.s. and in $L^1$ as }$N \rightarrow \infty$.$\hfill{\Box}$
\medskip

Note (see, e.g., \cite{G})  that $\kappa(p)$ is strictly positive and
finite for all $0<p<p_c$.
Furthermore, in \cite{Z}
the large deviations property of $K_N$ is formulated as follows
\medskip

\noindent
{\bf Theorem}
[\cite{Z}, Theorem 2] {\it
Given $\kappa(p)> \varepsilon>0$, there exist $\sigma_j(\varepsilon,
p)>0$
for $j=1,2$ such that
\[
\lim_{n\rightarrow \infty}\frac{-1}{|B(n)|}\log {\bf P}\left(
  \frac{K_N}{|B(n)|} \geq \kappa(p)+ \varepsilon \right) =
\sigma_1(\varepsilon,
p)
\]
and
\[
\lim_{n\rightarrow \infty}\frac{-1}{|B(n)|}\log {\bf P}\left(
  \frac{K_N}{|B(n)|} \leq \kappa(p)- \varepsilon \right) =
\sigma_2(\varepsilon,
p).
\]
}
\medskip

\noindent
This theorem  implies that
for all $0<\delta < \kappa(p)$ and all large $N$ there is a positive
constant $\sigma=\sigma (\delta,
p)$ such that 
\begin{equation}\label{T(N)} 
 {\bf P}\left\{ \left| \frac{K_N}{|B(N)|}- \kappa(p)
\right| > \delta \right\}
  \leq e^{ -\sigma \,|B(N)|}.
\end{equation}

Define for any $k \geq 1$ and $x\geq 0$ an indicator function:
\begin{displaymath}
    I_{k}(x) = \left\{ \begin{array}{ll}
            1, & \mbox{ if } x=k, \\
            0, & \mbox{ otherwise. }
        \end{array} \right.
\end{displaymath}

\begin{prop}\label{P2} For any fixed $k\geq 1$
\begin{equation}\label{mu} 
  \frac{1}{K_N} \sum _{i=1}^ {K_N}  I_{k}(|X_i|) \
  {\rightarrow} \ \frac{1}{\kappa(p)} \  \frac{{\bf P}
  \{|C|=k\}}{k} =:{ \mu}(k)
\end{equation}
a.s. and in $L^1$   as $N \rightarrow \infty$.
\end{prop}

\noindent
{\bf Proof.} Let $C(z)$, $z \in B(N)$, denote a connected
(in $G_N^{(s)}(p)$) cluster
which contains vertex $z$. Then we write
\begin{equation}\label{N1} 
\frac{1}{K_N} \sum _{i=1}^ {K_N}  I_{k}(|X_i|) =
\frac{|B(N)|}{K_N} \, \frac{1}{k} \ \frac{1}{|B(N)|}\sum _{z \in B(N)}  I_{k}(|C(z)|).
\end{equation}
By the ergodic theorem
\begin{equation}\label{conv} 
\frac{1}{|B(N)|}\sum _{z \in B(N)}  I_{k}(|C(z)|) \ {\rightarrow} \  {\bf P}
  \{|C|=k\}
\end{equation}
  $a.s.$ as $N \rightarrow \infty$. This in turn implies that
  convergence (\ref{conv}) holds  in $L^1$ as well,  since
  $$0 \leq \frac{1}{|B(N)|} \sum_{z\in B(N)} I_k(|C(z)|)
 \leq 1.$$
Hence, statement (\ref{mu}) follows by (\ref {N1}) and  (\ref{P1}).
$\hfill{\Box}$

Given a collection of clusters ${\bf X}$ defined in (\ref{t}),   we
introduce another graph $ {\widetilde G}_{N}({\bf X}, p,c)$ as follows.
The set of vertices of $ {\widetilde G}_{N}({\bf X},p,c)$ we denote
$\{v_1, \ldots , v_{K_N}\}$.  Each vertex $v_i$
 is
said to be of 
type $|X_i|$,
which means that 
 $v_i$ corresponds to the set of $|X_i|$ connected vertices in
 $B(N)$. We shall also call any vertex $v_i$ of $ {\widetilde G}_{N}({\bf X}, p,c)$
a {\it macro-vertex}, and  write sometimes 
\begin{equation}\label{v}
  v_i=X_i.
\end{equation}
With this notation the type of a macro-vertex $v_i$ is simply the cardinality of set
$v_i=X_i$.
The space of 
 the types of macro-vertices
 is  $S=\{1,2, \ldots  \}$. According to (\ref{mu}) the
 distribution of type of a (macro-)vertex in graph $ {\widetilde
   G}_{N}({\bf X},p,c)$
 converges to measure $\mu$ on $S$.
  The edges between the vertices of $ {\widetilde G}_{N}({\bf X},p,c)$ are
 presented independently
 with probabilities induced by the original graph
 $G_{N}(p,c)$. More precisely,
the probability of an edge between any two vertices $v_i$ and $v_j$
of types $x$
and $y$ correspondingly,  is 
 ${\widetilde p}_{xy}(N)$
 introduced in (\ref{p}).
Clearly, this construction provides a one-to-one correspondence
between the connected components in the graphs $ {\widetilde
  G}_{N}({\bf X},p,c)$ and $
  G_{N}(p,c)$: the number
  of the connected components is the same
  for
  both
  graphs, as well as the number of the involved vertices from $B(N)$
in two corresponding components. In fact, considering conditionally on
${\bf X}$ graph 
 $ {\widetilde
  G}_{N}({\bf X},p,c)$ we neglect only those long-range edges from 
$G_{N}(p,c)$, which connect vertices within each $v_i$, i.e., the
vertices which  are already
connected
 through the short-range edges.

  Consider now
\begin{equation}\label{p3}
    {\widetilde p}_{xy}( N) = 1-\left(1-\frac{c}{|B(N)|} \right)^{xy}
    =: \frac{\kappa ' _N (x,y)}{|B(N)|}.
\end{equation}
Observe that if $x(N)\rightarrow x$ and $y(N)\rightarrow y$ then
\begin{equation}\label{ka}
    \kappa ' _N (x(N), y(N))\rightarrow  c  xy
\end{equation}
 for all $x, y \in S$.
In order to place our
  model into the framework of the inhomogeneous random graphs
  from \cite{BJR} let us introduce another (random) kernel
  \[
    \kappa _{K_N} (x,y)=\frac{K_N}{|B(N)|}  \kappa_N '(x, y),
\]
so that we can rewrite the probability ${\widetilde p}_{xy}(N)$
  in a graph
  $ {\widetilde
  G}_{N}({\bf X},p,c)$ taking into account the size of the graph:
\begin{equation}\label{pij}
    {\widetilde p}_{xy}( N) = \frac{\kappa _{K_N} (x, y)}{K_N}.
\end{equation}
(We use notations from \cite{BJR} whenever it is appropriate.)
According to (\ref{P1}) and (\ref{ka}),
if $x(N)\rightarrow x$ and $y(N)\rightarrow y$ then
\begin{equation}\label{kap}
    \kappa _{K_N} (x(N), y(N))\stackrel{a.s.}{\rightarrow}   \kappa (x,
    y):= c\kappa(p)xy\ \ \
\end{equation}
as $N \rightarrow \infty$ for all $x, y \in S$ .

Hence,  
in view of Proposition \ref{P2}
we conclude that
conditionally on $K_N=t(N)$, where $t(N)/|B(N)| \rightarrow {\bf E}(|C|^{-1})$,
 our model 
falls into  the so-called "rank 1 case" of the general inhomogeneous
random graph model $G^{\cal V}(t(N),\kappa _{t(N)})$ with a vertex space
$${\cal V}=(S,\mu, (v_1, \ldots, v_{t(N)})_{N\geq 1})$$
(see \cite{BJR}, Chapter 16.4). 
Note, that according to (\ref{X}) function $\mu(k)$ (defined in (\ref{mu})) decays
exponentially, which implies
\begin{equation}\label{c1}
\kappa \in L^1(S \times S,  \mu\times \mu).
\end{equation}
 Furthermore, it is not difficult  to verify with a
help of (\ref{P1}) and Proposition \ref{P2}  that
for any $t(N)$ such that $t(N)/|B(N)| \rightarrow {\bf
  E}(|C|^{-1})$
\begin{equation}\label{cond}
\frac{1}{t(N)} {\bf E}\{ e({\widetilde G}_{N}({\bf X},p,c))|K_N=t(N)\}
\rightarrow
\frac{1}{2}\sum_{y=1}^{\infty} \sum_{x=1}^{\infty}\kappa (x,y) {\mu}(x)
{\mu}(y),
\end{equation}
where $e(G)$ denotes the number of edges in a graph $G$.
According to
Definition 2.7 from \cite{BJR},  under 
the conditions (\ref{cond}), (\ref{c1}) and (\ref{kap})
the sequence of  kernels $\kappa_{t(N)}$
(on the countable space $S \times S$)
is called {\it graphical} on
${\cal V}$
with limit $\kappa$.

\subsection{A branching process related to $ {\widetilde
  G}_{N}({\bf X},p,c)$.}
Here we closely follow the approach from \cite{BJR}.
We shall use a well-known technique of branching processes to reveal the
connected component in graph $ {\widetilde
  G}_{N}({\bf X},p,c)$. Recall first the usual algorithm of finding a
connected component. Conditionally on the set of macro-vertices, take any vertex $v_i$ to be the root. Find all
the vertices $\{ v^1_{i_1},v^1_{i_2},...,v^1_{i_n}\}$ connected to this vertex $v_i$ in the graph $ {\widetilde
  G}_{N}({\bf X},p,c)$, call them the first generation of $v_i$, and
then mark $v_i$ as "saturated". Then for each
non-saturated but already revealed vertex, we find all the vertices
connected to them but which have not been used previously. 
We continue this process until we end
up with a tree of saturated vertices.

Denote  $\tau_{N}(x)$ the set of the macro-vertices in
the tree constructed according to the  above algorithm with the
root at a vertex of type $x$.

It is plausible to think (and in our case it is correct, as will be
seen below) that
this algorithm with a high
probability as $N\rightarrow \infty$ reveals a tree of the offspring
of the
following multi-type 
 Galton-Watson process with type space $S=\{1,2, \ldots\}$: at any step,
 a particle of type $x \in S$ is replaced in the next
 generation by a set of particles where the number of  particles of
 type $y$ has a Poisson distribution $Po(\kappa
 (x,y) {\mu}(y))$.
 Let $\rho(x)$ denote the probability that a particle of type $x$
 produces an infinite population.

\begin{prop}\label{P4}  The function $\rho(x)$, $x \in S$, is the maximum
  solution to
\begin{equation}\label{rho}
\rho(x) = 1- e^{-\sum_{y=1}^{\infty} \kappa
 (x,y) \mu(y)\rho(y) }.
\end{equation}
  \end{prop}

\noindent  
{\bf Proof.}
We have
\begin{equation*}
\sum_{y=1}^{\infty}\kappa(x,y)\mu(y)= c{\bf E}(|C|^{-1})x \frac{1}{{\bf E}(|C|^{-1})} =cx< \infty \text{ for
any $x$},
\end{equation*}
which together with (\ref{c1})
verifies that the conditions of
Theorem 6.1 from \cite{BJR} are satisfied, and the result (\ref{rho})
follows by
this theorem. $\hfill{\Box}$

\bigskip

Notice that it also follows by the same Theorem 6.1 from \cite{BJR}
that $\rho(x)>0$ for all $x \in S$ if and only if
\begin{equation}\label{cr12}
  c {\bf E}(|C|^{-1}) \sum_{y=1}^{\infty}y^2\mu(y)=
  c \, {\bf E}(|C|^{-1}) \sum_{y=1}^{\infty}y^2
\frac{1}{{\bf
  E}(|C|^{-1})} \  \frac{{\bf P}
  \{|C|=y\}}{y}=c \ {\bf E}|C|>1;
\end{equation}
otherwise, $\rho(x)= 0$ for all $x \in S$. Hence, formula  (\ref{cr1})  for the
critical value $c^{cr}(p)$ follows from (\ref{cr12}).

As we showed above, conditionally on $K_N$ so that
$K_N/|B(N)|\rightarrow
{\bf E}(|C|^{-1})$, the
sequence $\kappa_{K_N}$ is graphical on $\cal
V$. Hence, the conditions of Theorem 3.1 from \cite{BJR} are satisfied
and we derive (first,
conditionally on $K_N$, and therefore unconditionally)
that 
\[
\frac{C_1({\widetilde G}_{N}({\bf X}, p,c))}{K_N} \stackrel{P}{\rightarrow} \rho,
\]
where $\rho=\sum_{x =1}^{\infty}\rho(x)\mu(x)$. This
together with (\ref{P1}) 
implies
\begin{equation}\label{MG}
\frac{C_1({\widetilde G}_{N}({\bf X}, p,c))}{|B(N)|}
\stackrel{P}{\rightarrow} {\bf E}(|C|^{-1}) \ \rho.
\end{equation}
Notice that here $C_1({\widetilde G}_{N}({\bf X}, p,c))$ is the number
of {\it macro-vertices} in the largest connected component of
${\widetilde G}_{N}({\bf X}, p,c)$.

\subsection{On the distribution of types of vertices in $ {\widetilde
  G}_{N}({\bf X},p,c)$.}

 Given a collection of clusters ${\bf X}$ (see (\ref{t}))
we define for all $1 \leq k \leq |B(N)|$
\[
 {\cal N} _k= {\cal N} _k ({\bf X})= \sum _{i=1}^ {K_N}  I_{k}(|X_i|). \]
In words,  ${\cal N} _k$ is  the number of (macro-)vertices of type $k$ in the set of
vertices of graph $ {\widetilde
  G}_{N}({\bf X},p,c)$. We shall prove here a useful result on the
distribution of ${\cal N}=({\cal N} _1, \ldots , {\cal N} _{K_N})$.

\begin{lem}\label{L1}Set 
$${\widetilde \mu}(k) =  \sum_{n = k}^{\infty}{\bf P} \{ |C|=n \}  $$
 and fix $\nu > 2$ arbitrarily. Then for any fixed $\varepsilon >0$
\begin{equation}\label{x} 
{\bf P}
  \Big\{  | {\cal N}_k/K_N - \mu  (k)| > \varepsilon \,k^{\nu}\, {\widetilde \mu}(k) \ \
  \mbox{ for some } 1\leq k \leq |B(N)| \Big\}
  =o(1)
\end{equation}
as $N \rightarrow \infty$.
\end{lem}
\noindent
{\bf Proof.} Let us fix $\varepsilon >0$ arbitrarily. Define a
constant $L_0$ so that $\varepsilon  L_0^{\nu}\,= {\bf E}(|C|^{-1})$.
Then for all $k>L_0$
\begin{equation}\label{10}
 \varepsilon \, k^{\nu}\, {\widetilde \mu}(k) > \mu
(k),
\end{equation}
and for any $L>L_0$
\begin{equation}\label{A1}
{\bf P}
  \{  | {\cal N}_k/K_N - \mu  (k)| > \varepsilon \,  k^{\nu}\, {\widetilde \mu}(k)  \ \
  \mbox{ for some } 1\leq k \leq  |B(N)|\}
\end{equation}
  \[ \leq {\bf P} \{ \max_{1\leq i \leq K_N } |X_i| >
  L \} 
\]
\[+ {\bf P}
  \{  | {\cal N}_k/K_N - \mu  (k)| > \varepsilon \, k^{\nu}\, {\widetilde \mu}(k)   \ \
  \mbox{ for some } 1\leq k \leq L \}.
\]
We shall choose later on an appropriate  $L=L(N)$ so that we will be able to
bound from above by $o(1)$ (as $N\rightarrow \infty$)
each of the summands on the right in (\ref{A1}).

First we derive
\begin{equation}\label{A2n}
{\bf P} \{  \max_{1\leq i \leq K_N} |X_i| > L \}
=
  {\bf P} \{ 
    \max_{z \in B(N) } |C(z)| > L \}
\leq
  |B(N)| \, {\bf P} \{ |C| > L\},
\end{equation}
where $C(z)$ is an open cluster containing $z$.
For a further reference we note here, that according to (\ref{X})
for any $0<\alpha < \zeta(p)$ there is constant $b>0$ such that
\begin{equation}\label{A5}
{\bf P} \{|C|\geq L \} \leq be^{-\alpha L}
\end{equation}
for all $L \geq 1$, which together with (\ref{A2n}) implies, in particular,
that
\begin{equation}\label{A3}
{\bf P} \{  \max_{1\leq i \leq K_N} |X_i| > \frac{2}{\zeta(p)}\log |B(N)|\} 
\rightarrow 0
\end{equation}
as $N \rightarrow \infty$.

Now we consider the last term in (\ref{A1}).
Let us define for any $0<\delta <{\bf E}(|C|^{-1})$  an event
\begin{equation}\label{calA}
{\cal A}_{\delta , N} =\left \{ \left| \frac{K_N }{|B(N)|}- {\bf E}(|C|^{-1})\right|
\leq
\delta \right\}.
\end{equation}
Recall that according to (\ref{T(N)})
\begin{equation}\label{AA}
{\bf P} ({\cal A}_{\delta , N}) \geq 1- e^{-\sigma |B(N)|}
 =1-o(1)
\end{equation}
as $N \rightarrow \infty$. Then we can bound
the last term in (\ref{A1}) as follows
\begin{multline}\label{A2}
{\bf P}
  \Big\{  | {\cal N}_k/K_N - \mu  (k)| > \varepsilon \, k^{\nu} {\widetilde \mu}(k) \ \
  \mbox{ for some } 1\leq k \leq L \Big\}
\\
\leq 
{\bf P}
  \Big\{  | {\cal N}_k/K_N - \mu  (k)| > \varepsilon \, k^{\nu} {\widetilde \mu}(k) \ \
  \mbox{ for some } 1\leq k \leq L_0 \Big\}
\\
+{\bf P}
 \Big\{   \Big( {\cal N}_k/K_N  > \varepsilon \,k^{\nu}\, {\widetilde \mu}(k) + \mu  (k)\ \
  \mbox{ for some } L_0 < k \leq L  \Big)\cap  {\cal A}_{\delta, N}
  \Big\} + {\bf P}\{\overline{{\cal A}_{\delta, N}}\} 
\\
\leq {\bf P}
 \Big\{  \Big(  {\cal N}_k/K_N  > \varepsilon \,k^{\nu}\, {\widetilde \mu}(k) + \mu  (k)\ \
  \mbox{ for some } L_0 < k \leq L  \Big)\cap  {\cal A}_{\delta, N}
  \Big\} +o(1),
\end{multline}
as $N\rightarrow \infty$ where the last inequality follows by Proposition \ref{P2} and bound (\ref{AA}).
Write
\begin{equation}\label{A7}
{P}(k):  =   {\bf P}
\bigg\{ \left(
 \frac{{\cal N}_k}{K_N}- \mu  (k)
  > \varepsilon \, k^{\nu}\, {\widetilde \mu}(k)\right)\cap {\cal A}_{\delta ,
    N} \bigg\}.
\end{equation}
Clearly, we have by  (\ref{A2}):
\begin{equation}\label{A7n}
{\bf P}
 \Big\{  
  | {\cal N}_k/K_N - \mu  (k)| > \varepsilon \, k^{\nu} {\widetilde \mu}(k) \ \
  \mbox{ for some } 1\leq k \leq L
  \Big\} \leq \sum_{k=L_0+1}^{L}{ P}(k) +o(1),
\end{equation}
as $N\rightarrow \infty$.
Substituting now (\ref{A7n}) and (\ref{A2n}) into (\ref{A1})
we derive
\begin{equation}\label{Be3}
{\bf P}
  \{  | {\cal N}_k/K_N - \mu  (k)| > \varepsilon \,  k^{\nu}\, {\widetilde \mu}(k)  \ \
  \mbox{ for some } 1\leq k \leq  |B(N)|\}
\end{equation}
  \[ \leq |B(N)| \, {\bf P} \{ |C| > L\} + \sum_{k=L_0+1}^{L}{ P}(k) +o(1)
 \leq |B(N)|{\widetilde \mu}(L)   \,  + \sum_{k=L_0+1}^{L}{ P}(k) +o(1)
\]
as $N\rightarrow \infty$.

Next we shall find an upper bound for ${ P}(k)$.
Due to the definition (\ref{calA}) of 
${\cal A}_{\delta , N}$, we have
\begin{equation}\label{compl}
{ P}(k)
 \leq {\bf P}
\bigg\{ {\cal N}_k> (\kappa (p)-\delta)|B(N)|
\Bigl(\varepsilon \, k^{\nu}\, {\widetilde \mu}(k) + \mu (k) \Bigr)
 \bigg\}.
\end{equation}
We shall use the following special case of the Talagrand's inequality.

\begin{prop}\label{PT} 
For every $r \in \mathbb{R}$ and $t \geq 0$
\begin{equation}\label{TI}
{\bf P} \biggl\{ {\cal N}_k \leq r-t \biggr\} {\bf P}
\biggl\{ {\cal N}_k  \geq r \biggr\} \leq \exp \left\{-\frac{t^2}{8dkr}\right\}.
\end{equation}
\end{prop}
\noindent
{\bf Proof.}
We shall derive this result as a corollary to the Talagrand's inequality 
 \cite{T} which we cite here from the book \cite{JLR}, p. 40.

\bigskip

\noindent
{\bf Theorem.} [Talagrand's Inequality]
{\it Suppose that $Z_1, ..., Z_n$ are independent random variables taking
their values in some sets $\Lambda_1, ... , \Lambda_n$,
respectively. Suppose further that $W=f(Z_1, ..., Z_n)$, where $f:
\Lambda_1 \times ... \times \Lambda_n \rightarrow \mathbb{R}$ is a
function such that, for some constants $c_k$, $k=1,...,n$, and some
function $\Psi$, the following two conditions hold:

\medskip

1) If $z, z' \in \Lambda = \prod_1^n \Lambda_i$ differ only in the
  $i-{th}$ coordinate, then $|f(z)-f(z')| \leq c_i$.

\medskip

2) If $z \in \Lambda$ and $r \in \mathbb{R}$ with $f(z)\geq r$,
  then there exists a set $J \subseteq \{1,...,n\}$ with $\sum_{i \in
    J}c_i^2 \leq \Psi (r)$, such that for all $y \in \Lambda$ with
  $y_i=z_i$ when $i \in J$, we have $f(y)\geq r$.

Then, for every $r \in \mathbb{R}$ and $t \geq 0$,
\begin{equation}\label{Tal}
{\bf P} (W \leq r-t) {\bf P} (W \geq r) \leq e^{-t^2/  4 \Psi(r)}.
\end{equation}
}
\medskip

We shall show now that
function ${\cal N}_k$ satisfies the conditions of this theorem.
Let $\{e_1, \ldots , e_n\}$ be the set of all edges from the lattice
${\bf Z}^d$ which have both end points 
 in $B(N)$.
Define 
\begin{displaymath}
  Z_i=  \left\{ \begin{array}{ll}
            1, & \mbox{ if $e_i$ is open in $G_N^{s}(p)$}, \\
            0, & \mbox{ if $e_i$ is closed in $G_N^{s}(p)$ .}
        \end{array} \right.
\end{displaymath}
According to the definition of our model, $Z_i \in Be(p)$, $i=1, \ldots , n$,
are independent random variables, and 
\[{\cal N}_k={\cal N}_k(Z_1, \ldots, Z_n)\] 
since ${\cal N}_k$ is the number of the components of size $k$
(open $k$-clusters) 
 in $G_N^{s}(p)$, which is defined completely by $Z_1, \ldots, Z_n$. 
Furthermore, it is clear that removing or adding just one edge in
$G_N^{s}(p)$,
may increase or decrease by at most one the number of $k$-clusters. 
Hence, the first condition of the Talagrand's inequality  is satisfied with $c_i=1$ for all $1\leq i\leq n$: 
if configurations $ z, z' \in \{0,1\}^n$ 
 differ only in the $i^{th}$ coordinate, then 
\[|{\cal N}_k(z)-{\cal N}_k(z')| \leq 1.\]

Next we check that the second condition is fulfilled as well, and we
shall determine the
function $\Psi$. Assume, $z \in \{0,1\}^n$ corresponds such configuration of the
edges in $B(N)$  that  
${\cal N}_k(z) \geq
r$, for some $r\geq 1$, i.e., 
 there are at least $r$ clusters of size $k$. 
Let $\{e_j, j\in J\} \subset \{e_1, \ldots, e_n\}$ be a set of 
edges which have at least one common vertex with a set
of exactly $r$
(arbitrarily chosen out of ${\cal N}_k(z) $) 
clusters of size $k$. Clearly, $|J|\leq 2dkr$, and
for any $z'\in \{0,1\}^n$
with $z'_j=z_j$ if $j\in J$, we have
\[{\cal N}_k(z') \geq
r,\]
proving that the second condition of the Talagrand's inequality is satisfied as well with $\Psi(r)=2dkr$, since
\begin{equation}\label{psir}
\sum_{i \in J}c_i^2 = |J| \leq 2dkr.
\end{equation}
Hence, the inequality (\ref{TI}) follows by (\ref{Tal}). \hfill$\Box$
\bigskip

Set now
\[k_N=(\kappa (p)-\delta)|B(N)|,\]
and consider  the inequality (\ref{TI})
with 
\begin{equation}\label{rt}
\begin{array}{rl}
r & = k_N
\Bigl( \varepsilon \, k^{\nu}\, {\widetilde \mu}(k) + \mu (k) \Bigr), \\ \\
r-t & =  k_N
\Bigl( \frac{\varepsilon}{2} \, k^{\nu}\, {\widetilde \mu}(k) + \mu
(k) \Bigr).
\end{array}
\end{equation}
First we derive for any $k>L_0$ (in which case $\varepsilon k^{\nu}\, {\widetilde \mu}(k) \geq \mu (k)$)
\begin{equation}\label{no1}
{\bf P} \biggl\{ {\cal N}_k
\leq r -t\biggr\}
\geq {\bf P} \biggl\{  {\cal N}_k\leq \frac{3}{2}
k_N \mu(k)\biggr\}
 \geq 1-\frac{{\bf E}  {\cal N}_k }{\frac{3}{2}k_N \mu (k)},
\end{equation}
where we used the Chebyshev's inequality.
Recall that by Proposition \ref{P2} (and (\ref{P1})) 
\[\frac{{\bf E}  {\cal N}_k }{|B(N)|} \rightarrow \kappa(p)\mu (k)\]
as $N\rightarrow \infty$. Hence, choosing $0< \delta \leq \kappa(p)/10$
we have 
\[\frac{{\bf E}  {\cal N}_k }{\frac{3}{2}k_N \mu (k)}=\frac{{\bf E}  {\cal N}_k }{
\frac{3}{2}(\kappa(p)-\delta)|B(N)| \mu (k)} \leq 3/4\]
for all large $N$, which together with (\ref{no1}) implies
\begin{equation}\label{no3}
 {\bf P} \biggl\{ {\cal N}_k
\leq r -t \biggr\} \geq \frac{1}{4}
\end{equation}
for all large $N$.
Using the last bound in the Talagrand's inequality (\ref{TI}) with $r$
and $t$ defined in (\ref{rt}), we derive for all large $N$ when $k>L_0$ (and therefore 
$\varepsilon k^{\nu}\, {\widetilde \mu}(k) \geq \mu (k)$)
\begin{equation}\label{TI1}
 {\bf P}
\biggl\{ {\cal N}_k  \geq r \biggr\} \leq 
\left({\bf P} \biggl\{ {\cal N}_k \leq r-t \biggr\} \right)^{-1}
\exp \left\{-\frac{t^2}{8dkr}\right\}
\end{equation}
\[
\leq 4 \exp \left\{-\frac{\left( 
\frac{\varepsilon}{2} k_N \, k^{\nu}\, {\widetilde \mu}(k)
\right)^2}{8dk \left( k_N
\Bigl( \varepsilon \, k^{\nu}\, {\widetilde \mu}(k) + \mu (k) \Bigr)\right)}\right\}
\]
\[
\leq 4 \exp \left\{-\frac{ 
\varepsilon k_N \, k^{\nu}\, {\widetilde \mu}(k)
}{64dk }\right\}
=4 \exp \left\{-\frac{ 
\varepsilon (\kappa(p)-\delta)
}{64d } |B(N)|k^{\nu -1}\, {\widetilde \mu}(k)
\right\}.
\]
Substituting (\ref{TI1}) 
into (\ref{compl}) we get
\[
{ P}(k)
 \leq 4 \exp \left\{-a |B(N)|k^{\nu -1}\, {\widetilde \mu}(k)
\right\},\]
where
\[a:= \frac{ 
\varepsilon (\kappa(p)-\delta)
}{64d } .\]
The last bound combined with (\ref{Be3}) yields
\begin{equation}\label{Be6}
{\bf P}
 \Big\{  
  | {\cal N}_k/K_N - \mu  (k)| > \varepsilon \, k^{\nu} {\widetilde \mu}(k) \ \
  \mbox{ for some } 1\leq k \leq |B(N)|
  \Big\} 
\end{equation}
\[
\leq |B(N)| \,  {\widetilde \mu}(L) +
4 \sum_{k=L_0}^{L} \exp \left\{-a |B(N)|k^{\nu -1}\, {\widetilde \mu}(k)
\right\}  
+o(1),\]
as $N\rightarrow \infty$ for any $L \geq L_0$.

Next we shall show that for any $\Delta >0$ one can choose a finite constant $L_0$ and numbers $L=L(N)$  such that
\begin{equation}\label{Be2}
 \sum_{k=L_0}^{L(N)}
\exp \left\{-a |B(N)|k^{\nu -1}\, {\widetilde \mu}(k)
\right\} < \Delta, 
\end{equation}
for all large $N$, and
\begin{equation}\label{Be4}
|B(N)| \,  {\widetilde \mu}(L(N)) \rightarrow 0, 
\end{equation}
as $N\rightarrow \infty$. This together  with 
 (\ref{Be6})  will clearly imply the statement of Lemma.

We claim that both (\ref{Be2}) and  (\ref{Be4}) hold with 
\begin{equation}\label{def}
L(N)=\min \left\{k : k^{\alpha} {\widetilde \mu}(k)<\frac{1}{|B(N)|} \right\},
\end{equation}
where 
\[\alpha = \frac{\nu-2}{2}\]
is positive by the assumption of Lemma.
Observe, that $k^{\alpha} {\widetilde \mu}(k) \rightarrow 0$,
as $k\rightarrow \infty$ for any fixed $\alpha$ due to the exponential decay
(\ref{X}). 
This yields that $L(N) \rightarrow \infty$ as $N\rightarrow \infty$, which in turn 
implies that there exists
\begin{equation*}
\lim_{N \rightarrow \infty} {\widetilde \mu}\bigl(L(N)\bigr)|B(N)| <
 \lim_{N \rightarrow \infty} L(N)^{-\alpha} =0,
\end{equation*}
and (\ref{Be4}) follows.

To prove (\ref{Be2})
first we note that by the definition (\ref{def}) of $L(N)$
\begin{equation}\label{bo}
(L(N)-1)^{\alpha}{\widetilde \mu} (L(N)-1) \geq \frac{1}{|B(N)|}.
\end{equation}
Recall that according to Lemma 6.102 from \cite{G} (p.139), for all $n,m\geq 0$
\begin{equation}\label{pin}
\frac{1}{m+n} {\bf P} (|C|=n+m) \geq p (1-p)^{-2}\frac{1}{m}{\bf P}
(|C|=m) \frac{1}{n} {\bf P} (|C|=n).
\end{equation}
When $m=1$ the inequality \eqref{pin} implies 
\begin{equation}
{\bf P} (|C|=n+1) \geq p (1-p)^{2(d-1)}{\bf P} (|C|=n),
\end{equation}
for all $n\geq 0$. This clearly yields
\begin{equation}\label{mun}
{\widetilde \mu} (L(N)) \geq p (1-p)^{2(d-1)} {\widetilde \mu} (L(N)-1).
\end{equation}
Notice that $ \gamma:= p (1-p)^{2(d-1)} \leq p <1$ for all $d\geq 1$.
Combining (\ref{bo}) with \eqref{mun} we immediately get
\begin{equation}\label{bo1}
L(N)^{\alpha} {\widetilde \mu} (L(N)) \geq  \frac{\gamma}{|B(N)|},
\end{equation}
and also by the definition (\ref{def}) for all $k<L(N)$
\begin{equation}\label{bo1*}
k^{\alpha} {\widetilde \mu} (k) \geq  \frac{1}{|B(N)|}
\geq  \frac{\gamma}{|B(N)|}.
\end{equation}
Making use of  (\ref{bo1}) and  (\ref{bo1*}) we derive
\begin{equation}\label{Be5}
 \sum_{k=L_0}^{L(N)}
\exp \left\{-a |B(N)|k^{\nu -1}\, {\widetilde \mu}(k)
\right\} 
\end{equation}
\[ \leq
\sum_{k=L_0}^{L(N)}\exp \left\{-a \gamma k^{\nu -1-\alpha
}\right\}
 \leq
a_1 \exp\left\{-a\gamma L_0^{\nu -2-\alpha} \right\},
\]
where $a_1$ is some positive constant independent of $L_0$. It is clear now, that for any 
$\Delta >0$ we can fix $L_0$ so large that (\ref{Be5}) would imply (\ref{Be2}), and
in the same time $L_0$ will satisfy (\ref{10}) and $L_0 <L(N)$.
This completes the proof of the lemma.
 \hfill$\Box$

\subsection{Proof of Theorem \ref{T} in the subcritical case $c  < c^{cr} (p)$  .}
Let us fix $0\leq p< p_c$ and then $c  < c^{cr} (p)$  arbitrarily.
Given ${\bf X}$ let again $v_i$ denote the macro-vertices with types
$|X_i|$, $i=1,2, \ldots,$
respectively, and let
${\widetilde
  L}$
denote a connected component in 
${\widetilde
  G}_{N}({\bf X},p,c)$. 
Consider now
for any positive constant $a$ and a function $w=w(N)$
\begin{equation}\label{J1}
{\bf P} \left\{  C_1\Big( G_N(pc)  \Big) > aw 
\right\}
=
{\bf P} \left\{ \max _{{\widetilde
  L}} \sum _{v_i \in {\widetilde
  L}}X_i
> a w \right\} .
\end{equation}
We know already from (\ref{MG})
that in the subcritical case
the
size (the number of macro-vertices) of any ${\widetilde
  L}$ is $o(N)$ with probability tending to one as $N\rightarrow \infty$.
Note that  when the kernel $\kappa(x,y)$
is not bounded uniformly in both arguments, which is our case,
it is not granted that the
largest component in the subcritical case is at most of order $\log
|B(N)|$ (see, e.g., discussion of 
Theorem 3.1 in \cite{BJR}). Therefore first we shall prove the following
intermediate result.

\begin{lem}\label{LS} If $c<c^{cr}(p$ then 
  \begin{equation}\label{S11}
{\bf P} \left\{  C_1\Big( {\widetilde
  G}_N({\bf X}, p,c)  \Big) > |B(N)|^{1/2}
\right\}
=o(1),
\end{equation}
as $N\rightarrow \infty$.
\end{lem}

\noindent
{\bf Proof.}
Let us fix $\varepsilon >0$ and $\delta >0$ arbitrarily and introduce
the following event
\begin{equation}\label{A19}
{\cal B}_N
={\cal A}_{\delta, N}
 \cap \left( \max _{1\leq i\leq K_N} |X_i| \leq \frac{2}{\zeta(p)}\log
   |B(N)| \right)
\end{equation}
 \[ \cap \left( \cap _{k=1}^{|B(N)|} \left\{\left|
       \frac{{\cal N}_k}{K_N} -\mu(k)\right|\leq \varepsilon k^{\nu}
     \widetilde 
     {\mu}(k)
\right\} \right).
\]
According to (\ref{AA}),  (\ref{A3}) and
(\ref{x}) we have
\begin{equation}\label{cb1}
{\bf P} \left\{ {\cal B}_N \right\} =1 -o(1)
\end{equation}
as $N \rightarrow \infty$.

Recall that  $\tau_{N}(x)$ denote the set of the macro-vertices in
the tree constructed according to the algorithm of revealing of
connected component described above.
Let $|\tau_{N}(x)|$ denote the
number of macro-vertices in
$\tau_{N}(x)$.
Then we easily derive
\begin{equation}\label{SA18}
 {\bf P} \left\{  C_1\Big( {\widetilde
  G}_N({\bf X}, p,c)  \Big) > |B(N)|^{1/2}
\right\}
\leq
{\bf P} \left\{ \max _{1\leq i\leq K_N} |\tau_{N}(|X_i|)| >  |B(N)|^{1/2}\mid {\cal B}_N 
\right\} +o(1)
\end{equation}
\[
\leq
|B(N)|\Big( \delta + {\bf E}(|C|^{-1}) \Big)
\sum_{k=1}^{|B(N)|} 
({\mu}(k)+ \varepsilon k^{\nu}\widetilde{\mu}(k) )
{\bf P} \left\{  
|\tau_{N}(k)| >  |B(N)|^{1/2}\mid {\cal B}_N 
\right\} +o(1)
\]
as $N \rightarrow \infty$.
We shall use the multi-type branching process introduced above
(Section 2.2)
to approximate the distribution of $|\tau_{N}(k)|$.
Let further ${\cal X}^{c,p}(k)$ denote the total number of the particles
(including the initial one) produced by
the branching process
starting with a single particle
of  type $k$. 
Observe that at each step of the exploration algorithm,
the number of new neighbours of $x$ of type $y$ has a
binomial distribution $ Bin(N_y',{\widetilde  p}_{xy}(N))$ where $N_y'$
is the number of remaining vertices of type $y$, so that 
$ N_y'
\leq {\cal N}_y $. 

We shall explore the following obvious 
relation between the Poisson and the binomial distributions. 
Let $Y_{n,p} \in Bin(n,p)$ and $Z_{a} \in Po(a)$, where $0<p<1/4$ and $a>0$.
Then for all $k\geq 0$
\begin{equation}\label{A21}
{\bf P}
  \{ Y_{n,p}=k\}
\leq (1+Cp^2 )^n \, 
{\bf P}\{ Z_{n\frac{p}{1-p}}=k
  \},
\end{equation}
where $C$ is some positive constant (independent of $n$, $k$ and $p$).
Notice that for all $1\leq x,y \leq \frac{2}{\zeta(p)}\log |B(N)|$ 
\begin{equation}\label{pB}
    {\widetilde p}_{xy}(N)= 1-\left(1-\frac{c}{|B(N)|} \right)^{x y}
    =\frac{c}{|B(N)|} \ {x y} \ (1+o(1)),
\end{equation}
and clearly,  ${\widetilde p}_{xy}(N) \leq 1/4$ for all large $N$.
 Therefore 
for any fixed 
positive $\varepsilon _1$
we can choose small $\varepsilon$ and $\delta$ in (\ref{A19}) so that
conditionally on ${\cal B}_N$ we have
\begin{equation}\label{pB1}
    N_y' \frac{{\widetilde p}_{xy}(N)}{1-{\widetilde p}_{xy}(N)}
 \leq  ({\mu} (y)+ y^{\nu} \varepsilon _1\widetilde{\mu} (y))\kappa (x,y)
\end{equation}
for all large $N$.

We shall use the following property of measure 
$$\mu(k)=\frac{1}{{\bf
  E}(|C|^{-1})} \  \frac{{\bf P}
  \{|C|=k\}}{k} $$
defined in Proposition \ref{P2}. Recall, that along with the 
result (\ref{X}) it is also proved in \cite{G}
that for all $0<p<p_c$
\begin{equation}\label{X1}
\zeta(p)=lim_{n \rightarrow \infty}\left(
-\frac{1}{n} \log {\bf P} \left\{|C|\geq n \right\}\right).
\end{equation}
Hence,
 (\ref{X}) and (\ref{X1}) immediately imply the existence and equality of the following limits for all $0<p<p_c$
\begin{equation}\label{X2}
\zeta(p)=lim_{n \rightarrow \infty}\left(
-\frac{1}{n} \log \mu(n)\right)
=lim_{n \rightarrow \infty}\left(
-\frac{1}{n} \log \widetilde{\mu}(n) \right),
\end{equation}
i.e., that both $\mu (n)$ and $\widetilde{\mu}(n)$ decay exponentially fast, and moreover with the same exponent
in the limit.
Let us write further
\[ \mu(y) =\mu _{p} (y), \ \ \ \ \ \widetilde{\mu} (y) =\widetilde{\mu}  _{p} (y), 
 \   \   \   \         \   \kappa (x,y)=   \kappa_{c,p} (x,y) =c\kappa(p)xy,\]
emphasizing dependence on $p$ and $c$. The result (\ref{X2}) allows us 
to choose  for any positive
$\varepsilon _2$ and 
$p<p'< p_c$ 
a positive $\varepsilon _1 = \varepsilon _1 (\varepsilon _2, p')$ such that
\begin{equation}\label{pB2}
   {\mu}_p (y)+ y^{\nu} \varepsilon _1
\widetilde{\mu}_p (y)
 \leq  
(1+\varepsilon _2)\mu_{p'}  (y).
\end{equation}
Setting now
$c':=(1+\varepsilon _2)\frac{\kappa(p) }{\kappa(p')}
 \, c$
we derive from  (\ref{pB1}) with a help of (\ref{pB2}), that conditionally on ${\cal B}_N$
with an appropriate choice of constants
\begin{equation}\label{pB3}
    N_y' \frac{{\widetilde p}_{xy}(N)}{1-{\widetilde p}_{xy}(N)}
\leq  
(1+\varepsilon _2)\mu_{p'}  (y)\kappa_{c,p}  (x,y)
=
\mu_{p'}  (y)\kappa _{c', p'} (x,y).
\end{equation}
Recall that above  we fixed $p$ and $c<c^{cr}(p)$, where $c^{cr}(p)$ 
is strictly decreasing and continuous in $p$. Furthermore, function 
$\kappa(p)$ is analytic on $[0,p_c)$.
Hence, we can choose
 $p'>p$ and $c':=(1+\varepsilon _2) \,  \frac{\kappa(p) }{\kappa(p')}c$
so that 
\begin{equation}\label{A22}
    c<c'<c^{cr}(p')<c^{cr}(p),
\end{equation}
and moreover $c'$ and $p'$ can be chosen arbitrarily close to $c$ and $p$, respectively.
Now according to (\ref{A21}) and (\ref{pB3}) 
\begin{equation}\label{n2}
{\bf P}\{ 
Y_{N_y' ,{\widetilde p}_{xy}(N)} \geq k
\}
\leq 
(
1+C{\widetilde p}_{xy}(N)^2 
)^{N_y' }\, 
{\bf P}
\{ 
Z_{N_y' 
\frac{{\widetilde p}_{xy}(N)}{1-{\widetilde p}_{xy}(N)}}
\geq k
  \} 
\end{equation}
\[
\leq (1+C{\widetilde p}_{xy}(N)^2 )^{|B(N)|}
{\bf P} \{ Z_{\mu_{p'}  (y)\kappa_{c', p'}  (x,y)}\geq k
  \}.
\]
Hence, if conditionally on ${\cal B}_N $
 at each (of at most $|B(N)|$) step
of the exploration algorithm which reveals $\tau_N(k)$, we replace the
$ Bin(N_y',{\widetilde  p}_{xy}(N))$ variable 
with the
$Po\left(\mu_{p'}  (y)\kappa_{c', p'}  (x,y)\right)$ one, we arrive at
the
following bound using branching process:
\begin{equation}\label{S12}
{\bf P}
\left\{ |\tau_N(k) | > |B(N)|^{1/2}\mid {\cal B}_N 
\right\}
\end{equation}
\[
\leq
\left(1+C \left(\max_{x,y \leq 2 \log |B(N)|/\zeta(p)} {\widetilde p}_{xy}(N)
\right)^2\right)^{|B(N)|^2}
 \, {\bf P} \left\{ {\cal X}^{c',p'} (k)> |B(N)|^{1/2}
\right\} .
\]
 This together with (\ref{pB}) implies
\begin{equation}\label{SA23}
{\bf P} \left\{ |\tau_N(k) | > |B(N)|^{1/2} \mid {\cal B}_N 
\right\} 
\leq
e^{b_1log |B(N)|)^4} \, 
{\bf P} \left\{
  {\cal X}^{c',p'} (k)> |B(N)|^{1/2}
\right\} ,
\end{equation}
 where
$b_1 $ is some positive constant. Substituting the last bound into (\ref{SA18})
we derive with a help of (\ref{pB2})
\begin{equation}\label{SA24}
{\bf P} \left\{   C_1\Big( {\widetilde
  G}_N({\bf X}, p,c)  \Big) >|B(N)|^{1/2}\right\} 
\end{equation}
\[
\leq
b_2|B(N)| e^{b_1(\log |B(N)|)^4}
\sum_{k=1}^{|B(N)|} 
 k \mu _{p'}(k)
{\bf P} \left\{ {\cal X}^{c',p'} (k)> |B(N)|^{1/2}
\right\}
 +o(1)
\]
as $N \rightarrow \infty$, where
$b_2$ is some positive constant. By the Markov's inequality
\begin{equation}\label{SMark}
{\bf P} \big\{{\cal X}^{c',p'} (k) > |B(N)|^{1/2}\big\} \leq z^{-|B(N)|^{1/2}}
{\bf
    E} z^{{\cal X}^{c',p'} (k)}
\end{equation}
for all $z\geq 1$.
Denote $h_z(k)={\bf E} z^{{\cal X}^{c',p'} (k)}$; then with a help of 
  (\ref{SMark}) we get from (\ref{SA24})
\begin{equation}\label{SA25}
{\bf P} \left\{  C_1\Big( {\widetilde
  G}_N({\bf X}, p,c)  \Big) > |B(N)|^{1/2}  \right\} 
\end{equation}
\[\leq b_2|B(N)| e^{b(\log |B(N)|)^4} z^{-|B(N)|^{1/2}}
\sum_{k=1}^{|B(N)|} 
 k \mu _{p'}(k)
h_z(k) + o(1).\]
Now we will show that there exists $z > 1$ such that the series 
\begin{equation}\label{Bz}
D_z (c', p')= \sum_{k=1}^{\infty} k \mu_{p'}(k) h_z(k)=\kappa(p')^{-1}
\sum_{k=1}^{\infty} {\bf P}_{p'}\{|C|=k\} h_z(k)
\end{equation}
 converge. This together with (\ref{SA25}) will clearly imply the
 statement of the lemma.

Note that  
function $h_z(k)$ (as a generating function for a  branching process) satisfies the following equation
\[
\begin{array}{ll}
h_z(k) & =z \exp {\left\{
\sum_{x=1}^{\infty} \kappa_{c', p'}  (k,x) \mu_{p'}(x) (h_z(x)-1 )
\right\}} \\
\\
& =z\exp \left\{
c'\kappa(p') k \left(
\sum_{x=1}^{\infty} x \mu_{p'}(x)h_z(x)-\kappa(p')^{-1}\right)\right\} \\
\\
& =z \exp \left\{ c'\kappa(p') k  
( B_z(c', p')-\kappa(p')^{-1})\right\}.
\end{array}
\]
Multiplying both sides by $k \mu_{p'} (k)$ and summing up over $k$ 
we find
\[
D_z(c', p') = \sum_{k=1}^{\infty} k \mu_{p'} (k)z \exp \left\{c'\kappa(p') 
 k ( D_z(c', p')-\kappa(p')^{-1})\right\} 
\]
\[= \sum_{k=1}^{\infty}  \kappa(p')^{-1} {\bf P}_{p'}\{|C|=k\}z \exp \left\{c'\kappa(p') 
 k ( D_z(c', p')-\kappa(p')^{-1})\right\},\]
where we also used the definition of $\mu(k)$.
Let us write  for simplicity $D_z=D_z(c, p)$.
Hence, as long as $D_z$ is finite, it  should satisfy equation
\begin{equation}\label{Az}
D_z = \kappa(p)^{-1} 
z\, {\bf E} e^{c |C| ( \kappa(p)D_z-1 )},
\end{equation}
which implies in turn that  $D_z$ is finite for some $z>1$ if and only if (\ref{Az}) has at least
one solution (for the same value of $z$). 
Notice that by the definition (\ref{Bz})
\begin{equation}\label{Sinit}
D_z \geq D_1 = \kappa(p)^{-1}=({\bf E}(|C|^{-1}))^{-1}
\end{equation}
for  $z \geq 1$. Let us fix $z> 1$ and consider equation
\begin{equation}\label{SAz1}
y/z  = \, {\bf E} e^{c|C|( y-1 )}=:F(y)
\end{equation}
for $ y \geq 1 $.
Using the property (\ref{X}) of the distribution of $|C|$ it is easy to derive 
that
function $F(y)$ is defined on $[0, \zeta(p)/c)$ where it is finite,
  increasing
  and has positive second derivative. Compute now
\begin{equation}\label{SAz2}
\frac{\partial }{\partial  y}F(y)|_{y=1 } = 
{c}{ {\bf E} |C| } =\frac{c}{c^{cr}}.
\end{equation}
Hence,  if $c<c^{cr}$ then there exists $z>1$ such that 
 there is a finite solution $y$ to (\ref{SAz1}), and therefore
(\ref{Az}) also has at least
one solution for some $z>1$. 
 Taking into account
condition (\ref{A22}),   we find that  $
D_z(c', p')$ is also finite for some $z>1$, which 
 finishes the proof of the lemma.
\hfill$\Box$

\bigskip

Now we are ready to complete the proof of (\ref{E1}), following almost the same arguments as in the proof of the previous lemma.
Let $S_{N}(x)= \sum_{X_i \in \tau_{N}(x)}|X_i|$ denote the
number of vertices from $B(N)$ which compose the macro-vertices of
$\tau_{N}(x)$.
Denote
\[ {\cal B}'_N := {\cal B}_N \cap
\left( C_1\Big( {\widetilde
  G}_N({\bf X}, p,c)  \Big)
\leq |B(N)|^{1/2}
\right).
\]
According to (\ref{cb1}) and Lemma \ref{LS}
we have
\[{\bf P} \left\{ {\cal B}'_N \right\}=1- o(1).\]
This allows us to  derive from (\ref{J1}) 
\begin{equation}\label{A18}
{\bf P} \left\{  C_1 \Big({
  G}_{N}(p,c) \Big) > a w \right\} \leq
{\bf P} \left\{ \max _{1\leq i\leq K_N} S_N(X_i) > a w \mid
  {\cal B}'_N 
\right\} +o(1)
\end{equation}
\[
\leq
|B(N)| \Big( \delta + {\bf E}(|C|^{-1})\Big)\sum_{k=1}^{|B(N)|} 
({\mu}(k)+ \varepsilon k^{\nu}\widetilde{\mu}(k)) 
{\bf P} \left\{  
S_N(k) > a w \mid {\cal B}'_N 
\right\} +o(1).
\]

Let now $S^{c,p}(y)$  denote the 
sum of types (including the one of the initial particle) in  the total progeny
of the introduced above branching process 
starting with initial particle of type $y$.
Repeating the same argument which led to (\ref{S12}),
we get the
following bound using the introduced branching process:
\[
{\bf P} \left\{ S_N(k) > a w \mid {\cal B}'_N 
\right\} 
\]
\[
\leq
\left(1+C \left(\max_{x,y \leq 2 \log |B(N)|/\zeta(p)} {\widetilde p}_{xy}(N)
\right)^2\right)^{b_1 |B(N)| \sqrt{|B(N)|}}
 \, {\bf P} \left\{ S^{c',p'} (k)> a w 
\right\} 
\]
as $N \rightarrow \infty$, where we take into account that we can
perform at most $\sqrt{|B(N)|}$ steps of exploration (the maximal possible
number of  macro-vertices in any ${\widetilde L}$ conditioned on ${\cal B}'_N$ ).
This  together with (\ref{pB}) implies
\begin{equation}\label{A23}
{\bf P} \left\{ \tau_N(k) > a w \mid {\cal B}'_N 
\right\} 
\leq
(1+o(1))
{\bf P} \left\{ S^{c',p'} (k)> a w 
\right\} 
\end{equation}
as $N\rightarrow \infty$. Substituting the last bound into (\ref{A18})
we derive
\begin{equation}\label{A24}
{\bf P} \left\{  C_1 \Big({
  G}_{N}(p,c) \Big) > a w \right\} 
\leq b
|B(N)| 
\sum_{k=1}^{|B(N)|} 
 k \mu _{p'}(k)
{\bf P} \left\{ S^{c',p'} (k)> a w
\right\}
 +o(1)
\end{equation}
as $N \rightarrow \infty$, where $b$ is some positive constant. 
Denote $g_z(k)={\bf E} z^{S^{c',p'} (k)}$; then 
similar to (\ref{SA25})
 we derive from (\ref{A24})
\begin{equation}\label{A25}
{\bf P} \left\{  C_1 \Big({
  G}_{N}(p,c) \Big) > a w(N) \right\} 
\leq
b
|B(N)| \sum_{k=1}^{|B(N)|} 
 k \mu _{p'}(k)
g_z(k) z^{-a w(N)} + o(1).
\end{equation}
We shall search for all  $z\geq 1$ for which the series 
$$A_z (c', p')= \sum_{k=1}^{\infty} k \mu_{p'}(k) g_z(k)=
\frac{1}{\kappa(p')}\sum_{k=1}^{\infty} {\bf P}_{p'}\{|C|=k\} g_z(k)$$ 
 converge.
Function $g_z(k)$ (as a generating function for a certain branching process) satisfies the following equation
\[
\begin{array}{ll}
g_z(k) & =z^k \exp {\left\{
\sum_{x=1}^{\infty} \kappa_{c', p'}  (k,x) \mu_{p'}(x) (g_z(x)-1 )
\right\}} \\
\\
& =z^k \exp \left\{
c'\kappa(p') k 
\sum_{x=1}^{\infty} x \mu_{p'}(x)(g_z(x)-1) \right\} \\
\\
& =z^k \exp \left\{ c'k(\kappa(p') 
A_z(c', p')-1)\right\}.
\end{array}
\]
Multiplying both sides by $k \mu_{p'} (k)$ and summing up over $k$ 
we find
\[
A_z(c', p') = \sum_{k=1}^{\infty} k \mu_{p'} (k)z^k \exp \left\{
 c'k(\kappa(p') 
A_z(c', p')-1)
\right\} .
\]
Denoting for simplicity $A_z=A_z(c, p) $, we can rewrite the last
equation as follows:
\begin{equation}\label{SAz}
A_z = \sum_{k=1}^{\infty} k \mu_{p} (k)z^k \exp \left\{
 ck(\kappa(p) A_z-1)\right\}.
\end{equation}
It follows from here (and the fact that $A_z\geq
1/\kappa(p)=A_1$ for all $z\geq 1$) that  if there exists $z>1$ for which the
series $A_z$ converge, it should satisfy by (\ref{X})
\begin{equation}\label{con1}
z<e^{\zeta(p)}.
\end{equation}
According to  (\ref{SAz}), as long as 
 $A_z$ is finite it  satisfies  the equation
\[
A_z = (\kappa(p) )^{-1}{\bf E}  \left( z^{|C|} \, e^{c|C|( \kappa(p) A_z-1 )}\right),
\]
which implies that $A_z$ is finite for some $z>1$
if and only if the last equation
has  at least one solution
\begin{equation}\label{init}
A_z \geq A_1 = 1/\kappa(p).
\end{equation}
Let us fix $z> 1$ and consider equation
\begin{equation}\label{Az1}
y = {\bf E} \left( z^{|C|} \, e^{c|C|( y-1)}\right)
\end{equation}
for $y>1$.
  It is easy to check that at least for some $y>1$ and $z>1$
function $$f(y,z):={\bf E} \left( z^{|C|} \, e^{c|C|( y-1)}\right)$$ is
  increasing, it has all the derivatives of the
second order, and $\frac{\partial ^2}{\partial y^2} f(y,z)>0$. Compute now
\begin{equation}\label{Az2}
\frac{\partial }{\partial  y}f(y,z)|_{y=1, z=1}= 
c {\bf E} |C|  =\frac{c}{c^{cr}}.
\end{equation}
Hence, if 
$c>c^{cr}$ there is no solution $y\geq 1$ to (\ref{Az1}) for any $z>1$.
On the other hand, if $c<c^{cr}$ then there exists $1<z_0<e^{\zeta(p)}$ such that for all
$1\leq z <z_0$ there is a finite solution $y\geq 1$ to (\ref{Az1}). 
We shall find $z_0$ as the (unique!) value for which function $y$ is tangent to $f(y,z_0)$ if $y\geq 1$. 

First we rewrite (\ref{Az1}) as follows. Set
\[a=\frac{1}{c}\log z, \]
then (\ref{Az1}) is equivalent to 
\begin{equation}\label{Az1*}
y = {\bf E} \, e^{c|C|( y-1+a)},
\end{equation}
which after the change $x=y-1+a$ becomes
\begin{equation}\label{Az1**}
x+1-a = {\bf E} \, e^{c|C|x}.
\end{equation}
Here on the right we have a convex function with a positive second
derivative (for all $x<\zeta(p)/c$).
Notice also that by the
assumption
\[
\frac{\partial }{\partial  x} {\bf E} \, e^{c|C|x} \mid_{x=0} = 
{\bf E} \,c |C| <1. \]
Hence, there exists unique $y_0>0$ such that
\begin{equation}\label{Az2*}
\frac{\partial }{\partial  x} {\bf E} \, e^{c|C|x} \mid_{x=y_0} = 
{\bf E} \,\Big( c |C|  e^{c|C| y_0} \Big)=1.
\end{equation}
Define now
\begin{equation}\label{a_0}
a_0 :=1+y_0 -{\bf E} \, e^{c|C|y_0},
\end{equation}
which is  strictly positive due to the preceeding argument. 
Clearly, function $x+1-a_0$ is tangent to ${\bf E} \, e^{c|C|x}$.
Hence, 
for all $a\leq a_0$ equation (\ref{Az1**}) has at least one solution,
which implies due to (\ref{Az1*})
that for all
\begin{equation}\label{z_0}
z\leq z_0 :=e^{c a_0}= \exp \{c(1+y_0 -{\bf E} \, e^{c|C|y_0})\}
\end{equation}
equation (\ref{Az1}) has also at least one finite solution $y>1$.
This yields in turn that $A_z$ is finite for all  $z\leq z_0$.

Now 
taking into account that $c'>c$ and $p'>p$ can be chosen arbitrarily close to $c$ and $p$, respectively, we derive from (\ref{A25})
that for all $1<z<z_0$ 
\begin{equation}\label{S1}
{\bf P} \left\{  C_1 \Big({
  G}_{N}(p,c) \Big) > a w(N) \right\} 
\leq b(z)
|B(N)|  z^{-a w(N)} + o(1)
\end{equation}
as $N\rightarrow \infty$, where $b(z)<\infty$. This implies that for any 
$$a>1/\log z_0= (c+cy_0 -{\bf E} \, ce^{c|C|y_0})^{-1}$$
and $w(N)=\log |B(N)|$
\begin{equation}\label{S2}
{\bf P} \left\{  C_1 \Big({
  G}_{N}(p,c) \Big) > a \log |B(N)| \right\} 
= o(1)
\end{equation}
as $N\rightarrow \infty$, which proves (\ref{E1}).
\hfill$\Box$

\bigskip

To conclude this section we comment
on the methods 
used here. It is shown in \cite{T3} that in the subcritical case of
the classical random graph model $G_{n, c/n}$  (i.e., $p=0$ in terms of our model) the same
 method of generating functions
leads to a constant which is exactly $\alpha (0,c)$
(see (\ref{aoc})). The last
 constant is known to be 
 the principal term for the asymptotics of the size of the
largest component (scaled to $\log N$) in the subcritical case.
This gives us hope that the constant
$\alpha(p,c)$
 is close to the optimal one also for $p>0$.

Similar methods were used in \cite {T2} for 
some class of inhomogeneous random graphs, and in \cite{BJR} for a general class of models.
Note, however, some difference with the approach in \cite{BJR}.
 It is assumed in \cite{BJR}, Section 12,
that the generating function
for the corresponding branching process with the initial state $k$ (e.g., our function $g_z(k)$, $k\geq 1$)
is bounded uniformly in $k$. 
As we prove here this condition is not always necessary: we need only convergence of the series $A_z$, while $g_z(k)$ is  unbounded in $k$ in our case. Furthermore, our approach allows one to construct constant $\alpha(p,c)$ as a function of the parameters of the model.

  \subsection{Proof of Theorem \ref{T} in the supercritical case.}

Let $\mathcal{C}_k$ denote the set of vertices in the
$k$-th largest component in graph ${G}_N(p,c)$, and conditionally on ${\bf X}$
let $\widetilde{\mathcal{C}}_k$ denote the set of macro-vertices in the
$k$-th largest component in graph $\widetilde{G}_N({\bf X}, p,c)$ (ordered in any way if there are ties).
Let also ${C}_k$ and $\widetilde{{C}}_k$  denote correspondingly, their sizes.
According to our construction  for any connected component $\widetilde{L}$ in  
 $\widetilde{G}_N({\bf X}, p,c)$ there is a unique component $L$ in
 ${G}_N(p,c)$ such that they are composed of the same vertices from
 $B(N)$, i.e., in the notations (\ref{v})
\[L=\cup _{X_i \in \widetilde{L}}\cup _{z \in X_i}\{z\}=:V(\widetilde{L}).\]
Next we prove that with a high probability the largest components in both graphs consist of the same vertices.

\begin{lem} \label{P5}
For any $0\leq p<p_c$ if $c>c^{cr}(p)$ then 
\begin{equation}\label{C}
{\bf P}\{{\mathcal{C}}_1 =V(\widetilde{\mathcal{C}}_1)\} =1-o(1)
\end{equation}
as $N \rightarrow \infty$.
\end{lem}
\noindent
{\bf Proof.} In a view of the argument preceeding this lemma we have
\[{\bf P}\{{\mathcal{C}}_1 \neq V(\widetilde{\mathcal{C}}_1)\} 
={\bf P}\{{\mathcal{C}}_1 = V(\widetilde{\mathcal{C}}_k) \mbox{ for some } k\geq 2\} .\]
According to Theorem 12.6 from \cite{BJR}, 
conditions of which are satisfied here, 
in the supercritical case conditionally on $K_N$ such that $K_N/|B(N)|
\rightarrow {\bf E}(|C|^{-1})$,
we have {\bf whp} 
$\widetilde{C}_2=O(\log(K_N))$, which by  Proposition
\ref{P1} implies $\widetilde{C}_2=O(\log |B(N)| )$ {\bf whp}.
Also we know already from (\ref{MG})
that in the supercritical case
  $\widetilde{C}_1=O(|B(N)|)$ {\bf whp}, and therefore
${C}_1=O(|B(N)|)$ {\bf whp}. Hence, 
for some positive constants $a$ and $b$
\begin{equation}\label{C1}
{\bf P}\{{\mathcal{C}}_1 \neq V(\widetilde{\mathcal{C}}_1)\} 
={\bf P}\{{\mathcal{C}}_1 = V(\widetilde{\mathcal{C}}_k) \mbox{ for some } k\geq 2\} 
\end{equation}
\[ \leq
{\bf P} \left\{\left( 
\max _{k\geq 2}|V(\widetilde{\mathcal{C}}_k)| >b  |B(N)|
\right) \cap \left( \max _{k\geq 2}\widetilde{{C}}_k < a \log |B(N)| \right)
\right\} 
+o(1).
\]
It follows from (\ref{A3})
that
\[{\bf P}\{\max _{1 \leq i \leq K_N} |X_i| \geq \sqrt{|B(N)|}\}=o(1)
\]
as $N \rightarrow \infty$.
Now we derive
\begin{equation}\label{C2}
{\bf P} \left\{\left( 
\max _{k\geq 2}|V(\widetilde{\mathcal{C}}_k)| >b |B(N)|
\right) \cap \left( \max _{k\geq 2}\widetilde{{C}}_k < a \log |B(N)|\right)
\right\} 
\end{equation}
\[ \leq 
{\bf P} \left\{\left( 
\max _{k\geq 2}|V(\widetilde{\mathcal{C}}_k)| >b |B(N)|
\right) \cap \left( \max _{k\geq 2}\widetilde{{C}}_k < a \log |B(N)|\right)
\cap \left( \max _{1 \leq i \leq K_N} |X_i| < \sqrt{|B(N)|}\right)
\right\} \]
\[ +o(1) \ \leq 
{\bf P} \left\{\sqrt{|B(N)|} \, a \log |B(N)| > b|B(N)|
\right\} +o(1)=o(1).\]
Substituting this bound into (\ref{C1}) we immediately get (\ref{C}). \hfill$\Box$

\bigskip

Conditionally on ${\mathcal{C}}_1 =V(\widetilde{\mathcal{C}}_1)$ we have
\begin{equation}\label{A26}
\begin{array}{ll}
\frac{C_1}{|B(N)|} & 
= \frac{1}{|B(N)|} \sum_{i=1}^{K_N} |X_i| {\bf 1}\{
X_i \in \widetilde{\mathcal{C}}_1 \} \\
\\
& = \frac{1}{|B(N)|} \sum_{i=1}^{K_N}\sum_{k=1}^{|B(N)|} k 
{\bf 1}
\{
|X_i|=k \}  {{\bf 1}}\{
X_i \in \widetilde{\mathcal{C}}_1 \} \\
\\
& = \frac{K_N}{|B(N)|} \sum_{k=1}^{|B(N)|} k \frac{1}{K_N}
\#\{
X_i \in \widetilde{\mathcal{C}}_1 : |X_i|=k\}. 
\end{array}
\end{equation}
Note that Theorem 9.10 from \cite{BJR} (together with (\ref{P1}) in our case)
 implies  that 
\begin{equation}\label{A32}
\nu_N(k):=\frac{1}{K_N} \# \{
X_i \in \widetilde{\mathcal{C}}_1 (N): |X_i|=k\} \stackrel{P}{\rightarrow}
 {
  \rho}( k) \mu (k)
 \end{equation}
for each $k \geq 1$, where ${
  \rho}( k)$ is the maximal solution to (\ref{rho}).

We shall prove below that also 
\begin{equation}\label{A27}
W_N:=\sum_{k=1}^{|B(N)|} k\nu_N(k)\stackrel{P}{\rightarrow}
 \sum_{k=1}^{\infty} k {
  \rho}(k) \mu (k) =:  \beta \, \left({\bf E}(|C|^{-1})\right)^{-1}.
\end{equation}
Observe that according to (\ref{rho}) constant
$\beta $ (defined above) is the maximal solution to
\[\beta = {\bf E}(|C|^{-1}) \sum_{k=1}^{\infty} k {
  \rho}(k) \mu (k)={\bf E}(|C|^{-1})\sum_{k=1}^{\infty} k 
\left(
1- e^{-\sum_{y=1}^{\infty} \kappa
 (k,y) \mu(y)\rho(y) }\right)
\mu (k)\]
\[= 1- {\bf E}\Big( e^{-c|C|\beta} \Big).
\]
 This proves that 
$\beta$  is the maximal root of (\ref{be}).
Then (\ref{A27}) together with (\ref{P1}) will allow us to
derive from 
(\ref{A26}) that for any positive $\varepsilon $
\[
{\bf P}
 \Big\{ \Big|\frac{ C_{1}\Big(G_N(p,c)\Big)}{|B(N)|} - \beta \Big| > \varepsilon \mid 
{\mathcal{C}}_1 =V(\widetilde{\mathcal{C}}_1) \Big\}
\rightarrow 0
\]
as $N \rightarrow \infty$. This combined with Lemma \ref{P5} would
immediately imply
\begin{equation}\label{C1/N}
\frac{ C_{1}\Big(G_N(p,c)\Big)}{|B(N)|} \stackrel{P}{\rightarrow}  \beta ,
\end{equation}
and hence the statement of the theorem follows.

Now we are left with proving (\ref{A27}). 
For any $1 \leq R < |B(N)|$ write $W_N:=W_N^{R}+w_N^{R}$, where
\[W_N^{R}:=\sum_{k=1}^{R} k\nu_N(k), \ \ 
w_N^{R}:=\sum_{k=R+1}^{|B(N)|} k\nu_N(k).
\]
By (\ref{A32}) we have for any fixed $R\geq 1$
\begin{equation}\label{A33}
W_N^{R} \ \stackrel{P}{\rightarrow} \ 
 \sum_{k=1}^{R} k {
  \rho}(k) \mu (k) 
\end{equation}
 as $N\rightarrow \infty$. Consider $w_N^{R}$. Note 
that for any $k\geq 1$
\begin{equation}\label{A34}
{\bf E}\nu _N(k)
\leq
{\bf E}
\frac{1}{K_N} \sum_{i=1}^{K_N} I_k(|X_i|)={\bf E}
 \frac{|B(N)|}{K_N} \frac{1}{k}
\frac{1}{|B(N)|}
 \sum_{z\in B(N)} I_k(|C(z)|)   ,
\end{equation}
where $C(z)$ denotes again a connected cluster in the bond percolation model on $B(N)$ with a probability $p$ of bound.
Using events ${\cal A}_{\delta, N}$ together with
bound (\ref{AA}),
we obtain from (\ref{A34}) for any fixed $0<\delta <{\bf E}(|C|^{-1})/2$
and $k\geq 1$
\[
{\bf E}\nu _N(k) \leq
 {\bf E} \left(
   \frac{|B(N)|}{K_N} {\bf 1}\{{\cal A}_{\delta, N}\}
\frac{1}{|B(N)|}
\sum_{z\in B(N)}I_k(|C(z)|)\right)
 +
{\bf E} \left( \frac{|B(N)|}{K_N} 
{\bf 1}\{\overline{{\cal A}_{\delta, N}}\}  \right)
\]
\[ \leq  \frac{1}{{\bf E}(|C|^{-1})-\delta } \
{\bf P}\{
 |C|=k \}+|B(N)|{\bf P}\{\overline{{\cal A}_{\delta, N}}\} . 
\] 
Bound (\ref{AA}) allows us to derive from here that
\begin{equation}\label{A35}
{\bf E}\nu _N(k) \leq A_1 ( {\bf P}\{
 |C|=k \}+e^{-a_1 |B(N)|})
\end{equation}
for some positive constants $A_1$ and $a_1$ independent of $k$ and $N$.
This together with (\ref{X}) yields
\begin{equation}\label{A36}
{\bf E} w_N^{R}=\sum_{k=R+1}^{|B(N)|} k 
{\bf E}\nu_N(k) \leq A_2 e^{-a_2 R}
\end{equation}
for some positive constants $A_2$ and $a_2$. 

Clearly, for any $\varepsilon >0$ we can choose $R_0$ so that for all
$R\geq R_0$
\[\sum_{k=R+1}^{\infty} k {
  \rho}(k) \mu (k) < \varepsilon /3,\]
and then we have
\begin{equation}\label{A37}
{\bf P} \{| W_N-\sum_{k=1}^{\infty} k {
  \rho}(k) \mu (k) |>\varepsilon\}
\end{equation}
\[={\bf P} \{|( W_N^R -\sum_{k=1}^{R} k {
  \rho}( k) \mu (k))
 +w_N^R-
\sum_{k=R+1}^{\infty} k {
  \rho}(k) \mu (k) |>\varepsilon\}\]
\[\leq {\bf P} \{| W_N^R -\sum_{k=1}^{R} k {
  \rho}(k) \mu (k)|>\varepsilon /3\}
 +{\bf P} \{w_N^R>\varepsilon /3\}. \]
Markov's inequality together with bound (\ref{A36}) gives us 
\begin{equation}\label{A38}
{\bf P} \{w_N^R>\varepsilon /3\} \leq \frac{3 {\bf E} w_N^{R}}{\varepsilon}
\leq \frac{3 A_2 e^{-a_2 R}}{\varepsilon}.
\end{equation}
Making use of (\ref{A38}) and (\ref{A33}) we immediately derive from 
(\ref{A37}) 
\begin{equation}\label{A39}
{\bf P} \{| W_N-\sum_{k=1}^{\infty} k {
  \rho}( k) \mu (k) |>\varepsilon\}
\leq o(1) + \frac{3 A_2 e^{-a_2 R}}{\varepsilon}
\end{equation}
as $N\rightarrow \infty$. Hence, for any given positive $\varepsilon$ and
$\varepsilon _0$ we can choose finite $R$ so large that  
\begin{equation}\label{A40}
\lim _{N \rightarrow \infty}{\bf P} \{| W_N-\sum_{k=1}^{\infty} k {
  \rho}(k) \mu (k) |>\varepsilon\}
< \varepsilon_0.
\end{equation}
This clearly proves statement (\ref{A27}), and therefore finishes the proof of 
the theorem.
\hfill$\Box$ 
\bigskip

\bigskip

\noindent
{\bf Acknowledgment} T.T.
 thanks MSRI for the hospitality at the beginning of this project.

\end{document}